\documentclass[a4paper,12pt]{amsart}

\usepackage{amsmath}
\usepackage{amssymb}
\usepackage{mathrsfs}
\usepackage{ifthen}
\usepackage{graphicx}
\usepackage[T1]{fontenc} %skandit

\nonstopmode \numberwithin{equation}{section}
\newtheorem{thm}[equation]{Theorem}
\newtheorem{cor}[equation]{Corollary}
\newtheorem{lem}[equation]{Lemma}
\newtheorem{prop}[equation]{Proposition}

\theoremstyle{definition}
\newtheorem{defn}{Definition}[section]

\newtheorem{prob}[equation]{Problem}

\newtheorem{innercustomthm}{Theorem}
\newenvironment{customthm}[1]
  {\renewcommand\theinnercustomthm{#1}\innercustomthm}
  {\endinnercustomthm}

\newcounter{minutes}
\divide\time by 60
\newcounter{hours}
\multiply\time by 60
\addtocounter{minutes}{-\time}
\newcounter {own}
\def\theown {\thesection       .\arabic{own}}

\newenvironment{pf}[1][]{%
 \vskip 3mm
 \noindent
 \ifthenelse{\equal{#1}{}}%
  {{\slshape Proof. }}%
  {{\slshape #1.} }%
 }%
{\qed\bigskip}

\newcounter{alphabet}
\newcounter{tmp}

\def\be{\begin{equation}}
\def\ee{\end{equation}}

\newcommand{\bee}{\begin{enumerate}}
\newcommand{\eee}{\end{enumerate}}

\newcommand{\blem}{\begin{lem}}
\newcommand{\elem}{\end{lem}}
\newcommand{\bthm}{\begin{thm}}
\newcommand{\ethm}{\end{thm}}
\newcommand{\bcor}{\begin{cor}}
\newcommand{\ecor}{\end{cor}}
\newcommand{\beg}{\begin{examp}}
\newcommand{\eeg}{\end{examp}}
\newcommand{\begs}{\begin{examples}}
\newcommand{\eegs}{\end{examples}}
\newcommand{\bdefe}{\begin{defin}}
\newcommand{\edefe}{\end{defin}}
\newcommand{\bprob}{\begin{prob}}
\newcommand{\eprob}{\end{prob}}
\newcommand{\bei}{\begin{itemize}}
\newcommand{\eei}{\end{itemize}}

\begin{document}

\title{Zalcman Conjecture for certain analytic and univalent functions}
%\title{Extreme Point Theory and its Applications to certain extremal problems}

\author{Vasudevarao Allu}
\address{Vasudevarao Allu,
Discipline of Mathematics,
School of Basic Sciences,
Indian Institute of Technology  Bhubaneswar,
Argul, Bhubaneswar, PIN-752050, Odisha (State),  India.}
\email{avrao@iitbbs.ac.in}

\author{Abhishek Pandey}
\address{Abhishek Pandey,
Discipline of Mathematics,
School of Basic Sciences,
Indian Institute of Technology  Bhubaneswar,
Argul, Bhubaneswar, PIN-752050, Odisha (State),  India.}
\email{ap57@iitbbs.ac.in}


\def\thefootnote{}
\footnotetext{ {\tiny File:~\jobname.tex,
printed: \number\year-\number\month-\number\day,
          \thehours.\ifnum\theminutes<10{0}\fi\theminutes }
} \makeatletter\def\thefootnote{\@arabic\c@footnote}\makeatother

\thanks{}

\maketitle
\pagestyle{myheadings}
\markboth{Vasudevarao  Allu and Abhishek Pandey}{ Zalcman Conjecture for certain analytic and univalent functions}

\begin{abstract}
Let $\mathcal{A}$ denote the class of analytic functions in the unit disk  $\mathbb{D}$  of the form $f(z)= z+\sum_{n=2}^{\infty}a_n z^n$
and $\mathcal{S}$ denote the class of functions $f\in\mathcal{A}$ which are  univalent ({\it i.e.}, one-to-one). In 1960s,   L. Zalcman conjectured that 
$|a_n^2-a_{2n-1}|\le (n-1)^2$ for $n\ge 2$, which 
implies the famous Bieberbach conjecture $|a_n|\le n$ for $n\ge 2$.
For $f\in \mathcal{S}$,   Ma \cite{Ma-1999} proposed a generalized Zalcman conjecture  
$$|a_{n}a_{m}-a_{n+m-1}|\le (n-1)(m-1)
$$
for $n\ge 2, m\ge 2$. Let $\mathcal{U}$ be the class of functions  $f\in\mathcal{A}$ satisfying 
$$
\left|f'(z)\left(\frac{z}{f(z)}\right)^2-1 \right|< 1 \quad\mbox{ for } z\in\mathbb{D}.
$$
and $\mathcal{F}$ denote the class of functions $f\in \mathcal{A}$   satisfying ${\rm Re\,}(1-z)^{2}f'(z)>0$ in $\mathbb{D}$.
%consider the class of close-to-convex functions with argument $0$ and with respect to Koebe function $k(z)=z/(1-z)^{2}$. That is,
%      $$\mathcal{F}=\{f\in \mathcal{A}: {\rm Re\,}(1-z)^{2}f'(z)>0 z\in \mathbb{D}\}.$$
In the present paper, we prove the  Zalcman conjecture and generalized Zalcman conjecture   for the class $\mathcal{U}$ using extreme point theory.
We aslo  prove  the  Zalcman conjecture and generalized Zalcman conjecture  for the class $\mathcal{F}$ for the initial coefficients.

%We consider the functional $|\mu a_n^2-a_{2n-1}|$, for real $\mu$ and determine 
%the sharp upper bound of this functional for functions in the closed convex hull of starlike and convex functions of order $\alpha$, $0\le\alpha<1$ 
%by using extreme point theory.

\end{abstract}

\section{Introduction and Preliminaries}\label{Introduction}

Let $\mathcal{H}$ denote the class of analytic functions in the unit disk $\mathbb{D}:=\{z\in\mathbb{C}:\, |z|<1\}$. Here  $\mathcal{H}$ is 
a locally convex topological vector space endowed with the topology of uniform convergence over compact subsets of $\mathbb{D}$. Let $\mathcal{A}$ denote the class of functions 
$f\in \mathcal{H}$ such that $f(0)=0$ and $f'(0)=1$. Let $\mathcal{S}$ 
denote the class of functions $f\in\mathcal{A}$ which are univalent ({\it i.e.}, one-to-one) in $\mathbb{D}$. 
If $f\in\mathcal{S}$ then $f(z)$ has the following representation
\begin{equation}\label{p4_i001}
f(z)= z+\sum_{n=2}^{\infty}a_n z^n.
\end{equation}
A function $f\in\mathcal{A}$ is called starlike (convex respectively) if $f(\mathbb{D})$ is starlike with respect to the origin (convex respectively). Let $\mathcal{S}^*$ and $\mathcal{C}$ denote the class of starlike and convex functions in $\mathcal{S}$ respectively. It is well-known that a function $f\in\mathcal{A}$ is in $\mathcal{S}^*$ if,  and only if, ${\rm Re\,}\left(zf'(z)/f(z)\right)>0$ for $z\in\mathbb{D}$. Similarly, a function $f\in\mathcal{A}$ is in $\mathcal{C}$ if,  and only if, ${\rm Re\,}\left(1+zf''(z)/f'(z)\right)>0$ for $z\in\mathbb{D}$. From the above it is easy to see that $f\in\mathcal{C}$ if, and only if, $zf'\in\mathcal{S}^*$.
 Given $\alpha\in(-\pi/2,\pi/2)$ and $g\in \mathcal{S}^*$, a function $f\in \mathcal{A}$ is said to be close-to-convex with argument $\alpha$ and with respect to $g$ if 
$$ {\rm Re\,} \left(e^{i\alpha}\frac{zf'(z)}{g(z)}\right)>0, \,\,\, z\in \mathbb{D}.$$
Let $\mathcal{K}_{\alpha}(g)$ denote the class of all such functions, and $$\mathcal{K}(g):=\bigcup\limits_{\alpha\in(-\pi/2,\pi/2)}\,\mathcal{K}_{\alpha}(g) \mbox{ and } \mathcal{K}_{\alpha}:= \bigcup\limits_{g\in\mathcal{S}^*}\,\mathcal{K}_{\alpha}(g)$$
be the classes of close-to-convex functions with respect to $g$, and close-to-convex function with argument $\alpha$, respectively. The class
$$\mathcal{K}:=\bigcup\limits_{\alpha\in(-\pi/2,\pi/2)}\,\mathcal{K}_{\alpha}=\bigcup\limits_{g\in \mathcal{S}^{*}}\, \mathcal{K}(g)$$
is the class of all close-to-convex functions. It is well-known that every close-to-convex function is univalent in $\mathbb{D}$. Geometrically, $f\in \mathcal{K}$ means that the complement of the image domain $f(\mathbb{D})$ is the union of non-intersecting half- lines. These standard classes are related by the proper inclusions $\mathcal{C}\subsetneq\mathcal{S^*} \subsetneq\mathcal{K}\subsetneq\mathcal{S}$. \\

For $0<\lambda\le1$, let $\mathcal{U}(\lambda)$ be the class of  functions $f\in\mathcal{A}$ satisfying 
$$\left |f'(z)\left(\frac{z}{f(z)}\right)^{2}-1\right|<\lambda  \quad\mbox {  for }  z\in \mathbb{D}.
$$
Since $f'(z)(z/f(z))^{2} \ne 0$ in $\mathbb{D}$,  it follows that every function in the class $\mathcal{U}(\lambda)$ is non-vanishing in $\mathbb{D}\setminus \{0\}$. We set $\mathcal{U}:=\mathcal{U}(1)$.  
It is known that functions  in $\mathcal{U}(\lambda)$ are locally univalent and functions in the class $\mathcal{U}$ are univalent (see \cite{Aksentev-1958}).
Furthermore, Aksentev \cite{ABU-Muhana-1992} and Ozaki and Nunokawa \cite{ozaki-1972}  have shown that  functions in $ \mathcal{U}(\lambda)$ are univalent, {\it i.e.},
$\mathcal{U}(\lambda) \subseteq \mathcal{S}$  for  $0<\lambda\le 1$.\\

It is worth to mention here some important definitions and results related to the subsets of $\mathcal{H}$.
\begin{defn}
A set $\mathcal{F}\subseteq \mathcal{H}$ is normal if each sequence $\{f_{n}\}$ in $\mathcal{F}$ has a subsequence $\{f_{n_{k}}\}$ which converges to a function $f \in \mathcal{H}$ uniformly on every compact subset of $\mathbb{D}$.
\end{defn}
\begin{defn}
A set $\mathcal{F}\subseteq\mathcal{H}$ is locally bounded if for each point $a\in\mathbb{D}$ there are constants $M$ and $r>0$ such that for all $f \in\mathcal{H}$, 
$$|f(z)|\le M \mbox{ for } |z-a|<r.$$
That is, $\mathcal{F}$ is locally bounded if,  about each point $a\in\mathbb{D}$ there is a disk on which $\mathcal{F}$ is uniformly bounded.
\end{defn} 
\begin{lem}\cite[Lemma 2.8, p. 153]{Conway-1973}
A set $\mathcal{F} \subseteq\mathcal{H}$ is locally bounded if, and only if, for each compact set $K\subset\mathbb{D}$ there is a constant $M$ such that $$|f(z)|\le M$$ 
for all $f\in\mathcal{F}$ and $z\in K$. 
\end{lem}

\begin{thm} \cite[Montel's Theorem 2.9, p. 153]{Conway-1973} \label{Abhi-Vasu-thm}
A family $\mathcal{F} \subseteq\mathcal{H}$ is normal if, and only if, is locally bounded
\end{thm}
\begin{cor} \cite[Corollary 2.10, p. 154]{Conway-1973} \label{Abhi-Vasu-cor}
A set $\mathcal{F} \subseteq \mathcal{H}$ is compact if, and only if, it is closed and locally bounded.
\end{cor}

\begin{thm}\cite[Theorem 2.6 (Growth Theorem)]{Duren-book-1983}\label{Abhi-Vasu-P1-Growth}
For each $f\in \mathcal{S}$,
$$\frac{r}{(1+r)^{2}}\le |f(z)|\le \frac{r}{(1-r)^{2}},\,\,\,\,   |z|=r\,<\,1.$$
For each $z\in \mathbb{D}$, $z\ne0$, equality occurs if, and only if, $f$ is a suitable rotation of Koebe function.
\end{thm}
For $f$, $g \in \mathcal{H}$, we say that $f$ is subordinate to $g$, written as $f \prec g$ or $f(z) \prec g(z)$, if there 
exists an analytic function $\omega: \mathbb{D} \rightarrow \mathbb{D}$ with $\omega(0)=0$ such that 
$f(z)= g(\omega(z))$ for $z\in\mathbb{D}$. Furthermore, if $g$ is univalent in $\mathbb{D}$ then $f \prec g$ if, and only if, $f(0)= g(0)$ and $f(\mathbb{D})\subseteq g(\mathbb{D})$.
 %We say that $f$ is majorized by $g$ in $\mathbb{D}$ if $|f(z)|\le|g(z)|$ for each $z\in\mathbb{D}$. In other words, $f$ is majorized by $g$ in $\mathbb{D}$ if there exists an analytic function $\omega: \mathbb{D} \rightarrow \overline{\mathbb{D}}$ such that $f=\omega g$.
 If $\mathscr{G}\subseteq\mathcal{H}$, we use the notation $s(\mathscr{G})=\{f: f\prec g \mbox{ for some } g\in\mathscr{G}\}$. 
%and $m(\mathscr{G})=\{f: f \mbox{ is majorized by $g$ for some } g\in\mathscr{G}\}$.
 If $\mathscr{G}$ is a compact subset of $\mathcal{H}$ then it is 
not difficult to show that $s(\mathscr{G})$ is compact subset of $\mathcal{H}$ (for 
instance, see \cite[Lemma 5.19]{Hallenbeck-MacGregor-1984}).\\

Suppose $X$ is a linear topological vector space and $V \subseteq X$.  A point $x\in V$ is called an extreme point of $V$ if it has no representation of the form $x=ty+(1-t)z, 0<t<1$ 
as a proper convex combination of two distinct points $y$, $z\in V$. We denote $EV$ the set of extreme points of $V$. The convex hull of a set $V \subseteq X$ is the smallest convex set containing $V$. The closed convex 
hull denoted by ${\overline{co}}V$ is defined as the intersection of all closed convex sets containing $V$.  
That is,  the closed convex hull of $V$ is the smallest closed convex set containing $V$, which  is the closure of the convex hull of $V$. The Krein-Milman 
Theorem asserts that every compact subset of a locally convex topological space is contained in the closed convex hull of its extreme 
points (see, for instance, \cite{Dunford-1958}). For a general reference and for many important results on this topic, we refer to \cite{Hallenbeck-MacGregor-1984}. \\

As a first step for application of the knowledge of extreme point of these 
classes Brickman {\it et al.} \cite{Brickman-Macgregor-Wilken-1971} pointed out the following general results.

\begin{customthm}{A}\label{p4-theorem001}
Let $\mathscr{G}$ be a compact subset of $\mathcal{H}$ and $J$ be a complex-valued continuous linear functional on $\mathcal{H}$. Then
$\max\{ {\rm Re\, }J(f): f\in \overline{co}\, \mathscr{G}\} = \max\{ {\rm Re\, }J(f): f\in \mathscr{G}\} = \max\{ {\rm Re\, }J(f): f\in E\overline{co}\, \mathscr{G}\}$.
\end{customthm}

\begin{defn}\label{Abhi-Vasu-P1-Defn}
If $\mathcal{F}$ is a convex subset of $\mathcal{H}$ and $J:\mathcal{H}\to\mathbb{R}$ then $J$ is called convex on $\mathcal{F}$ provided 
that $J(tf+(1-t)g)\le tJ(f)+(1-t)J(g)$ whenever $f, g\in \mathcal{F}$ and $0\le t\le1$.
\end{defn} 

\begin{customthm}{B}\label{p4-theorem005}
Let $\mathscr{G}$ be a compact subset of $\mathcal{H}$ and $J$ be a real-valued, continuous and convex functional 
on $\overline{co}\, \mathscr{G}$. Then $\max\{ J(f): f\in \overline{co}\, \mathscr{G}\} = \max\{ J(f): f\in \mathscr{G}\} = \max\{ J(f): f\in E\overline{co}\, \mathscr{G}\}$.
\end{customthm}

The proof of these two results can be found in \cite[Theorem 4.5, Theorem 4.6]{Hallenbeck-MacGregor-1984}. In order to solve such linear 
extremal problems over $\mathscr{G}$, it suffices to solve them over the smaller class $E\overline{co}\, \mathscr{G}$. This reduction 
thereby becomes an effective technique for solving various linear extremal problems. Using this technique we solve the Zalcman conjecture for
the class $\mathcal{U}$.\\

% The coefficient bounds for the families $s(\mathcal{C})$ and $s(\mathcal{K})$ were first obtained by Rogosinski \cite{Rogosinski-1943} and Robertson \cite{Robertson-1965} 
%respectively. For the investigation of coefficient bounds for the families 
%$s(\mathcal{S^*})$, $m(\mathcal{S^*})$, $s(\mathcal{C})$, $m(\mathcal{C})$, $s(\mathcal{K})$ and $m(\mathcal{K})$ 
%with arguments using extreme point methods, we refer to \cite{Hallenbeck-MacGregor-1974,MacGregor-1967,MacGregor-1972} 
%(see also \cite{Hallenbeck-MacGregor-1984}). 

%Next we deal with a nonlinear extremal problem. 
In 1960s, L. Zalcman posed a conjecture that if a function $f\in\mathcal{S}$ is given by (\ref{p4_i001}) then
\begin{equation}\label{p4_i005}
|a_n^2-a_{2n-1}|\le (n-1)^2 \quad\mbox{ for } n\ge2,
\end{equation}
the equality holds only for the Koebe function $k(z)=z/(1-z)^2$ or its rotation.
It is important to note that the remarkable Zalcman conjecture implies the celebrated Bieberbach conjecture $|a_n|\le n$ for $f\in\mathcal{S}$  
(see \cite{Brown-Tsao-1986}). A well-known consequence of the area theorem  shows that  (\ref{p4_i005}) holds good for 
$n=2$ (see \cite{Duren-book-1983}). For $f\in\mathcal{S}$, Krushkal has proved the Zalcman conjecture for $n=3$ (see \cite{Krushkal-1995}) and recently 
for $n=4,5,6$ (see \cite{Krushkal-2010}). For a simple and elegant proof of Zalcman conjecture for the case $n=3$, we refer to \cite{Krushkal-2010}. 
The Zalcman conjecture for functions in the class $\mathcal{S}$ is still  open for $n>6$. However,  using complex geometry and universal Teichm\"{u}ller spaces Krushkal 
has  proved it for all $n\ge 2$ in his unpublished work \cite{Krushkal-Unpublished}.\\
 
The Zalcman conjecture has been proved affirmatively for certain special subclasses of $\mathcal{S}$, such as starlike functions, typically real functions, 
close-to-convex functions \cite{Brown-Tsao-1986,Ma-1988} and an observation also demonstrates that the Zalcman conjecture is asymptotically 
true (see \cite{Efraimidis-Vukotic-preprint}). Recently, Abu Muhana {\it et al.} \cite{Muhanna-Li-Ponnusamy-2015} solved Zalcman conjecture for the class 
$\mathcal{F}$ consists of  the family of  functions $f\in\mathcal{A}$ satisfying the condition
${\rm Re\, }(1+zf''(z)/f'(z))> -1/2$ for $z\in\mathbb{D}$. Functions in the class $\mathcal{F}$  are known to be convex in some direction (and hence close-to-convex and univalent) 
in $\mathbb{D}$.  In 1986, Brown and Tsao \cite{Brown-Tsao-1986} proved the Zalcman conjecture for the starlike functions and typically real functions. In 1988, Ma \cite{Ma-1988} 
proved that the Zalcman conjecture  for close-to-convex functions. For basic properties of starlike functions, typically real functions and close-to-convex functions 
we refer to \cite{Duren-book-1983, Vasu-book}.\\ %An observation also demonstrates that the Zalcman conjecture is asymptotically true (see \cite{Efraimidis-Vukotic-preprint}).\\

%%%%
%   General versions of Zalcman conjecture
%
%More general versions of Zalcman conjecture have also been considered by Brown and Tsao \cite{Brown-Tsao-1986}, and Ma \cite{Ma-1999} where the authors considered the functional $\Phi_{\lambda}(f)= \lambda a_n^2-a_{2n-1}$ for real parameter $\lambda$. As pointed out by Pfluger \cite{Pfluger-1976}, the  
% % % %quantity $\lambda a_n^2-a_{2n-1}$ arises quite naturally when we consider coefficient problems. Indeed, if $f\in\mathcal{S}$ then the coefficients of $(f(z^2))^{1/2}$ and $1/f(1/z)$ contain expressions of the form $\Phi_{\lambda}(f)$. Note that for the case $n=2$, the functional $\Phi_{\lambda}(f)$ reduces to  the classical Fekete-Szeg\H{o} functional $\Lambda_{\lambda}(f)= \lambda a_2^2- a_3$ which itself is a classical problem   having  long  and rich history in the literature (see \cite{Fekete-Szego-1933,Ali-Vasudevarao-2015,Koepf-1987,London-1993} and the references there in). Although the functional $\Phi_{\lambda}$ is not linear, the technique of extreme point theory can be used to find the maximum value of $\Phi_{\lambda}(f)$ over certain class $\mathscr{G}$. \\

In 1999, Ma \cite{Ma-1999} proposed a generalized Zalcman conjecture for $f\in \mathcal{S}$ that for $n\ge 2, m\ge 2$, $$|a_{n}a_{m}-a_{n+m-1}|\le (n-1)(m-1),$$
which is still an open problem. Ma  \cite{Ma-1999} has proved this generalized Zalcman conjecture for classes $\mathcal{S}^{*}$ and $\mathcal{S}_{\mathbb{R}}$.  Here $\mathcal{S}_{\mathbb{R}}$  denote 
the class of 
all functions in $\mathcal{S}$ with real coefficients. In 2017, Ravichandran and Verma \cite{Ravichandran - 2017} proved it for the classes 
of starlike and convex functions of given order and for the class of functions with bounded turning.\\

In the present paper, we prove the  Zalcman conjecture and generalized Zalcman conjecture   for the class $\mathcal{U}$ using extreme point theory.
We also  prove  the  Zalcman conjecture and generalized Zalcman conjecture  for the class $\mathcal{F}$ for the initial coefficients.
The organization of the paper is as follows. In Section \ref{Section-2} we prove that the class $\mathcal{U}(\lambda)$ for $0<\lambda\leq 1$ is compact.
In particular the class $\mathcal{U}$ is compact. In Section \ref{Section-3},  we will characterize the  closed convex hull of the class $\mathcal{U}$ and its extreme points. 
Then by using extreme point theory,  we prove the  Zalcman conjecture in Section \ref{Section-3} and 
generalized Zalcman conjecture    for the class $\mathcal{U}$ in Section \ref{Section-4}. We  prove  the  Zalcman conjecture and generalized Zalcman conjecture 
 for the class $\mathcal{F}$ for the initial coefficients in Section \ref{Section-5}.\\

Before we prove our main results we recall some important results which will play vital role in  our proofs. In 2016, Obradovi\'{c} {\it et al.} \cite{Obradvoic-2016} prove the following interested result.

\begin{prop} \cite{Obradvoic-2016}  \label{Abhi-Vasu-P1-Prop-001} 
 If $f\in \mathcal{U}(\lambda)$ for $0<\lambda\le 1$, then  for $z\in \mathbb{D}$, 
\begin{equation*}\label{Abhi-Vasu-P1-Equ-001abc} 
 \frac{f(z)}{z}\prec \frac{1}{(1-z)(1-\lambda z)}.
   \end{equation*}
\end{prop}

Let 
$$\mathcal{R}:=\left\{F\in \mathcal{H}:\overline{co}s(F)=\left\{\int_{|x|=1}F(xz)\,d\mu(x):\mu\in \wedge\right\}\right\},
$$ 
where $\wedge$ denote the set of probability measure on $\partial\mathbb{D}$.  We recall the following well-known result of Hallenbeck {\it et. al}  \cite{Hallenbeck-1989}. 

\begin{lem}\cite{Hallenbeck-1989}\label{Abhi-Vasu-P1-Lemma-Hallenbeck-1989} 
$\displaystyle  \frac{1}{(1-z)^{\alpha+i\beta}} \in \mathcal{R}$ if, and only if,  $\alpha\ge 1$ and $\beta=0$.
\end{lem}

Let $\mathcal{P}$ denote the class of all analytic functions $p$ in $\mathbb{D}$ with $p(0)=1$ satisfying ${\rm Re\,}p(z)>0$ in $\mathbb{D}$. 
Functions in the class $\mathcal{P}$ are called the $Carath\acute{e}odory$   functions and  can be expressed as 
\begin{equation}\label{Abhi-Vasu-P1-E}
p(z)=1+\sum\limits_{n=1}^{\infty}c_{n}z^{n}.
\end{equation}
\begin{lem}\cite{Hallenbeck-MacGregor-1984}\label{Abhi-Vasu-P1_Lemma-Hallenbeck-1984}
 $p\in \mathcal{P}$ if, and only if, there is a probability measure $\mu$ on $\partial\mathbb{D}$ such that  $$p(z)=\int_{|x|=1}\frac{1+xz}{1-xz}\, d\mu(x).$$ 
\end{lem} 
Equivalently, in view of the Lemma $\ref{Abhi-Vasu-P1_Lemma-Hallenbeck-1984}$, for $p\in \mathcal{P}$ given by $(\ref{Abhi-Vasu-P1-E})$ can be written as 
\begin{equation}\label{Abhi-Vasu-eq}
p(z)=1+\sum\limits_{n=1}^{\infty}c_{n}z^{n}=\int\limits_{0}^{2\pi}\frac{1+e^{it}z}{1-e^{it}z}\,\, d\nu(t).
\end{equation}
On comparing both the sides of $(\ref{Abhi-Vasu-eq})$ we obtain
\begin{equation}\label{Abhi-Vasu-P1-Equ}
c_{n}=2\int\limits_{0}^{2\pi}e^{int}\, d\nu(t).
\end{equation}

\begin{lem}\cite[Lemma 2.3, p. 507]{Ravichandran - 2015}\label{Lemma-Ravichandran-2015}
If $p(z)=1+\sum\limits_{k=1}^{\infty}c_{k}z^{k} \in \mathcal{P}$, then for all $n, m\in \mathbb{N}$,
\[|\lambda c_{n}c_{m}-c_{n+m}|\le
		         \begin{cases}
		             2,  & 0\le\lambda\le1\\[2mm]
		             2|2\lambda-1|,  & \mbox{ elsewhere }.\\
		         \end{cases}
		              \]
If $0<\lambda<1$, the inequality is	sharp for the function $p(z)=(1+z^{n+m})/(1-z^{n+m})$. In other cases, the inequality is sharp for the function $p(z)=(1+z)/(1-z)$.	              
\end{lem}

\section{Compactness of the set $\mathcal{U}(\lambda)$ }\label{Section-2}
\begin{thm}\label{Compctness}
For $0<\lambda\le1$, the class $\mathcal{U}(\lambda)$ is compact.
\end{thm}
\begin{pf}
In the view of Corollary $\ref{Abhi-Vasu-cor}$, it is enough to show
that $\mathcal{U}(\lambda)$  is closed and locally bounded. Let $\{f_{n}\}$ be a sequence in $\mathcal{U}(\lambda)$ which converges to 
$f$ uniformly on every compact subset of $\mathbb{D}$.  Clearly $f(0)=0$, $f'(0)=1$ and $f_{n}' \to f'$ uniformly on every compact subset of $\mathbb{D}$. 
Let
$$g_{n}(z)=\left(\frac{z}{f_{n}(z)}\right)^{2}f_{n}'(z)-1 \quad\mbox { and }\,  g(z)=\left(\frac{z}{f(z)}\right)^{2}f'(z)-1.
$$We aim to show that $g_{n} \to g$ uniformly on every compact subset of $\mathbb{D}$. Let $$h_{n}(z)=\frac{f_{n}(z)}{z} \quad\mbox { and }\, h(z)=\frac{f(z)}{z}.$$ Then  $h_{n}(z)\ne0$ and $h(z)\ne0$, for $ z\in \mathbb{D}$ and $ n\in \mathbb{N}$. 
Now we  prove that $h_{n} \to h$ uniformly on every compact subset of $\mathbb{D}$. 
To show this, it is enough to show that $h_{n}\to h$ uniformly on $\overline{D}=\{z\in\mathbb{C}:|z|\le r\}$, $0<r<1$,  where $D=\{z:|z|< r\}$.  Let
$$M_{n}=\sup_{z\in \overline{D}}|h_{n}(z)-h(z)|=\sup_{z\in \overline{D}}\frac{|f_{n}(z)-f(z)|}{|z|}.$$ Since $h_{n}(z)-h(z)$ is analytic function in $D$ and continuous on $\overline{D}$, so by maximum modulus theorem
 $$\max\{|h_{n}(z)-h(z)|:z\in \overline{D} \}=\max\{|h_{n}(z)-h(z)|:z \in \partial D\}.$$ %Hence $$M_{n}=\frac{M_{n}'}{r}, \mbox{ where } M_{n}'=\sup_{z\in \overline{D}}|f_{n}(z)-f(z)|$$
 Since $f_{n}$ converges to $f$ uniformly on $\overline{D}$ therefore, $M_{n}$ converges to $0$ as $n$ tends to $\infty$. %which gives us that $M_{n}$ converges to $0$. 
 Hence $h_{n}$ converges to $h$ uniformly on $\overline{D}$. Therefore, $h_{n}$ converges to $h$ uniformly on every compact subset of $\mathbb{D}$. Since $h_{n} \to h$ uniformly on every compact subset of $\mathbb{D}$ and $h_{n}(z)\ne0$, $h(z)\ne0$, for $ z\in \mathbb{D}$ and $ n\in \mathbb{N}$, it is not difficult to show that $1/h_{n}$ converges to $1/h$ uniformly on every compact subset of $\mathbb{D}$.\\ 
 
 If a sequence of continuous functions $f_{n}$ converge uniformly to a continuous function $f$ on some compact set, then $f_{n}$ is uniformly bounded on that compact set. In view of the above discussion,
  it is easy to see that
  $(1/h_{n})^2$ converges to $(1/h)^2$ uniformly on every compact subset of $\mathbb{D}.$
Also, $f_{n}'$ converges to  $f'$ uniformly on every compact subset of $\mathbb{D}$. Hence 
  $(1/h_{n})^{2}f_{n}' $ converges to $(1/h)^{2}f' $ unifromaly on every compact subset of $ \mathbb{D}$.
 Therefore $g_{n}$ converges to $g$ uniformly on every compact subset of $\mathbb{D}$.\\

Since  $|g_n(z)|<\lambda$ for each $n\in\mathbb{N}$, we prove that   $|g(z)|\le\lambda$. 
Suppose not, then there exists $z_{0} \in \mathbb{D}$ such that $|g(z_{0})|>\lambda$.
Let $\epsilon=|g(z_{0})|-\lambda$.  Then there exists $N\in \mathbb{N}$ such that 
$$ |g_{n}(z)-g(z)|<\epsilon=|g(z_{0})|-\lambda 
$$ 
for $ n\ge N$ and  $z\in \mathbb{D}$.  In particular, 
$$|g_{N}(z_{0})-g(z_{0})|<|g(z_{0})|-\lambda.
$$ 
Therefore,
$$|g(z_{0})|-|g_{n}(z_{0})|\le   |g_{N}(z_{0})-g(z_{0})|<|g(z_{0})|-\lambda,
$$ 
shows that  $|g_{n}(z_{0})|>\lambda$ which is a contradiction and hence  $|g(z)|\le\lambda$ in $\mathbb{D}$. 
If there exists some point $z_{0}\in \mathbb{D}$ such that $|g(z_{0})|=\lambda$ then by Maximum Modulus theorem, 
 $g$ must be a constant function, which is a contradiction. Therefore, $|g(z)|<\lambda$ for $z\in \mathbb{D}$ and hence $f\in \mathcal{U}(\lambda)$. This shows that  $\mathcal{U}(\lambda)$ is closed. 
In view of Theorem $\ref{Abhi-Vasu-P1-Growth}$, it is easy to observe that the class $\mathcal{S}$  is locally bounded. Since the class $\mathcal{U}(\lambda)\subseteq \mathcal{S}$ for $0<\lambda\leq 1$,  it follows that the class   $\mathcal{U}(\lambda)$ is also locally bounded. 
Thus, $\mathcal{U}(\lambda)$ is compact. 
\end{pf}

In particular,  for $\lambda=1$,  the class $\mathcal{U}$ is compact.

\section{Zalcman Conjecture for the class $\mathcal{U}$}\label{Section-3}

\begin{thm}\label{Close-to-convex-hull--Class-U}
$\overline{co}\, \mathcal{U}$ consists of all functions represented by  
\begin{equation*}
f(z)=\int_{|x|=1}\frac{z}{(1-xz)^{2}}\, d\mu(x), 
\end{equation*}
where $\mu \in \wedge$. Here $\wedge$ denotes the set of probability measure on $\partial{\mathbb{D}}$.
Further,  $E\overline{co}\, \mathcal{U}$ consists functions of the form 
\begin{equation*}
f(z)=\frac{z}{(1-xz)^{2}},\quad |x|=1.
\end{equation*}
\end{thm}
\begin{pf}
%[{\bf Proof of Theorem   \ref{Close-to-convex-hull--Class-U}}] 
%\begin{pf}
Let $f\in \mathcal{U}$ then in view of Proposition $\ref{Abhi-Vasu-P1-Prop-001}$,  we have 
$$\frac{f(z)}{z}\prec\frac{1}{(1-z)^{2}}.
$$ 
Let $F(z)=1/(1-z)^2$ then from Lemma $\ref{Abhi-Vasu-P1-Lemma-Hallenbeck-1989}$, we obtain
$$
\overline{co}s(F)=\left\{\int_{|x|=1}F(xz)\, d\mu(x):\mu\in \wedge\right\} \quad\mbox{ and } \quad E\overline{co}s(F)=\left\{F(xz):|x|=1 \right\}.
$$ 
Let 
$$\mathcal{G}:=\left\{\int_{|x|=1}\frac{z}{(1-xz)^{2}}\, d\mu(x):\mu\in\wedge \right\}
$$ 
then our aim is to  prove that  $\overline{co}\, \mathcal{U}=\mathcal{G}$. 
To prove this, we first prove that $\mathcal{G}$ is convex and compact. The fact that $\mathcal{G}$ is convex follows from the convexity of the set of probability measure on 
$\partial\mathbb{D}$ {\it i.e.},  convexity of $\wedge$.  To prove that $\mathcal{G}$ is compact, in the view of Corrolary $\ref{Abhi-Vasu-cor}$, we will show that $\mathcal{G}$ is  closed and locally bounded. 
 The fact that $\mathcal{G}$ is closed follows from the weak-star compactness of the set of probability measure on $\partial\mathbb{D}$.
To see that $\mathcal{G}$ is locally bounded,  let $|z|=r<1$ and $f\in \mathcal{G}$, then 
 $$f(z)=\int_{|x|=1}\frac{z}{(1-xz)^{2}}\, d\mu(x)
 .$$ 
 Since $1-|xz|\le |1-xz|$, we have $(1-r)^{2}\le|1-xz|^{2}$ for $|x|=1$ and hence
 \begin{equation}\label{Abhi-Vasu}
 \frac{|z|}{|1-xz|^{2}}\le \frac{r}{(1-r)^{2}}.
 \end{equation}
In view of $(\ref{Abhi-Vasu})$,  we obtain
$$|f(z)|=\left|\int_{|x|=1}\frac{z}{(1-xz)^{2}}\, d\mu(x)\right|\le \frac{r}{(1-r)^{2}}.
$$ 
Therefore $\mathcal{G}$ is locally bounded.\\

If $f\in \mathcal{U}$ then $f(z)/z\in s(F)$ which  implies that $f(z)/z\in \overline{co}s(F)$. That is,  
$$\frac{f(z)}{z}=\int_{|x|=1}\frac{1}{(1-xz)^{2}}\,\, d\mu(x)
$$ 
for some $\mu \in \wedge$ and hence,  
$$f(z)=\int_{|x|=1}\frac{z}{(1-xz)^{2}}\,\, d\mu(x)
$$ 
for some $\mu \in \wedge$. Therefore  $f \in \mathcal{G}$ and hence   $\mathcal{U}\subseteq \mathcal{G}$.
 Since $\mathcal{G}$ is closed and convex and $\overline{co}\, \mathcal{U}$ is the smallest closed convex set containing $\mathcal{U}$, it follows that
 $\overline{co}\, \mathcal{U}\subseteq\mathcal{G}$.  Since 
$$E\mathcal{G}=\left\{\frac{z}{(1-xz)^{2}}:|x|=1 \right\}
$$ 
and for each $x$ such that $|x|=1$, functions of the form $z/(1-xz)^2$ belong to $\mathcal{U}$, it follows that   $E\mathcal{G}\subseteq\mathcal{U}$.  Since $\mathcal{G}$ is compact and convex, we conclude that 
$\overline{co}\, \mathcal{G}=\mathcal{G}$  and hence by Krein-Milman Theorem,  it follows that 
$$\mathcal{G}\subseteq\overline{co}\, E\mathcal{G}\subseteq\overline{co}\, \mathcal{U}.
$$ 
Therefore,  $\mathcal{G}\subseteq\overline{co}\, \mathcal{U}$ and hence $\overline{co}\, \mathcal{U}=\mathcal{G}$ and 
$$E\overline{co}\, \mathcal{U}=\left\{\frac{z}{(1-xz)^{2}}:|x|=1\right\}.
$$
This completes the proof.
\end{pf}

Let $f\in \mathcal{S}$ be given by $(\ref{p4_i001})$. Then for fixed $n\in \mathbb{N}$, define the functional $\phi:\mathcal{S}\to \mathbb{C}$ defined by, $\phi(f)=a_{n}^{2}-a_{2n-1}$. 
 The rotations of $f\in\mathcal{S}$ be given by 
$$g(z)=e^{-i\theta}f(e^{i\theta}z)=z+\sum\limits_{n=2}^{\infty}A_{n}z^{n},
$$ 
where $A_{n}=a_{n}e^{i(n-1)\theta}$.  Since $\mathcal{S}$ is rotationally invariant, $g\in\mathcal{S}$. A simple computation shows that 
$$\phi(g)=A_{n}^{2}-A_{2n-1}=a_{n}^{2}e^{2i(n-1)\theta}-a_{2n-1}e^{2i(n-1)\theta}=e^{2i(n-1)\theta}\phi(f).
$$ 
This shows that $|\phi(g)|=|\phi(f)|$. That is,  $|a_{n}^{2}-a_{2n-1}|$ is invariant under rotations. Since $\mathcal{U}$ is rotationally invariant, maximizing $|a_{n}^{2}-a_{2n-1}|$ over $\mathcal{U}$ is equivalent to 
maximizing ${\rm Re\,}(a_{n}^{2}-a_{2n-1})$ over $\mathcal{U}$. It is easy to show that 
\begin{eqnarray}\label{Abhi-Vasu-array}
{\rm Re\,}(a_{n}^{2}-a_{2n-1}) & = & {\rm Re\,}(a_{n}^{2})-{\rm Re\,}(a_{2n-1})\\ \nonumber
&=&({\rm Re\,}(a_{n}))^{2}-({\rm Im\,}(a_{n}))^{2}-{\rm Re\,}(a_{2n-1})\\ \nonumber
& \le & ({\rm Re\,}(a_{n}))^{2}-{\rm Re\,}(a_{2n-1}).
\end{eqnarray}
 In view of $(\ref{Abhi-Vasu-array})$, we maximize $({\rm Re\,}(a_{n}))^{2}-{\rm Re\,}(a_{2n-1})$ over $\mathcal{U}$ to prove the Zalcman conjecture for the class $\mathcal{U}$.

\begin{thm}\label{Zalcman-Conjecture-Class-U}
Let  $f\in \mathcal{U}$  be given by $(\ref{p4_i001})$. Then  $|a_{n}^{2}-a_{2n-1}|\le (n-1)^{2}$ for $n \ge 2$. 
This inequality is sharp with equality for the Koebe function and its rotations {\it i.e.}, functions of the form $f(z)=z/(1-xz)^{2}$ where $|x|=1$.
\end{thm}
\begin{pf}
%[{\bf Proof of Theorem   \ref{Zalcman-Conjecture-Class-U}}] 
Since $\mathcal{U}\subseteq\mathcal{S}$, for the case $n=2$, the proof of Zalcman Conjecture holds good (see \cite{Duren-book-1983}) as consequence of the area theorem.
For the case  $n=3$, the Zalcman Conjecture has been proved  by Krushkal \cite{Krushkal-2010} for the class  $\mathcal{S}$. Therefore, it suffices to prove the Zalcman Conjecture for $n\geq 4$ for the class $\mathcal{U}$.  
For this, for fixed $n\in \mathbb{N}$ we define  the functional $J:\mathcal{S}\to \mathbb{R}$  by 
$$J(f)=({\rm Re\,}(a_{n}))^{2}-{\rm Re\,}(a_{2n-1}).
$$
We first prove that  $J$ is  convex on $\overline{co}\, \mathcal{U}$. Let $f,g\in\overline{co}\, \mathcal{U}$ be given by 
$f(z)=z+\sum_{n=2}^{\infty}a_{n}z^{n}$ and $g(z)=z+\sum_{n=2}^{\infty}b_{n}z^{n}$. In view of Definition $\ref{Abhi-Vasu-P1-Defn}$,   we show that  
$J(tf+(1-t)g)\le tJ(f)+(1-t)J(g)$.  Let
 $tf(z)+(1-t)g(z)=z+\sum_{n=2}^{\infty}A_{n}z^{n} 
 $,
 where $A_{n}=ta_{n}+(1-t)b_{n}$. A  computation shows that
\begin{eqnarray*}J(tf+(1-t)g)
&=&({\rm Re\,}(A_{n}))^{2}-{\rm Re\,}(A_{2n-1})\\ 
&=&({\rm Re\,}(ta_{n}+(1-t)b_{n}))^{2}-{\rm Re\,}(ta_{2n-1}+(1-t)b_{2n-1})\\ 
&=&({\rm Re\,}(ta_{n}+(1-t)b_{n}))^{2}-t{\rm Re\,}(a_{2n-1})-(1-t){\rm Re\,}(b_{2n-1})\\ 
&=& t^{2}({\rm Re\,}(a_{n}))^{2}+(1-t)^{2}({\rm Re\,}(b_{n}))^{2}+2t(1-t){\rm Re\,}(a_{n}){\rm Re\,}(b_{n})\\
&&\hspace*{0.25cm}  -t({\rm Re\,}(a_{n}))^{2}-(1-t)({\rm Re\,}(b_{n}))^{2}+t({\rm Re\,}(a_{n}))^{2}\\
&&\hspace*{0.25cm} +(1-t)({\rm Re\,}(b_{n}))^{2}-t{\rm Re\,}(a_{2n-1})-(1-t){\rm Re\,}(b_{2n-1})\\
&=& t^{2}({\rm Re\,}(a_{n}))^{2}-t({\rm Re\,}(a_{n}))^{2}(1-t)^{2}({\rm Re\,}(b_{n}))^{2}-(1-t)({\rm Re\,}(b_{n}))^{2}\\
&&\hspace*{0.25cm}+2t(1-t){\rm Re\,}(a_{n}){\rm Re\,}(b_{n}) t\left(({\rm Re\,}(a_{n}))^{2}-{\rm Re\,}(a_{2n-1})\right)\\
&&\hspace*{0.25cm} +(1-t)\left(({\rm Re\,}(b_{n})^{2}-{\rm Re\,}(b_{2n-1})\right)\\
&=& t(t-1)({\rm Re\,}(a_{n}))^{2}+t(t-1)({\rm Re\,}(b_{n}))^{2}-2t(t-1){\rm Re\,}(a_{n}){\rm Re\,}(b_{n})\\
&&\hspace*{0.25cm} +tJ(f)+(1-t)J(g) \\
&=& t(t-1)({\rm Re\,}(a_{n})-{\rm Re\,}(b_{n}))^{2}+tJ(f)+(1-t)J(g) \\
&=& tJ(f)+(1-t)J(g)-t(1-t)({\rm Re\,}(a_{n})-{\rm Re\,}(b_{n}))^{2}\\
&\le &  tJ(f)+(1-t)J(g).
\end{eqnarray*}
and hence $J$ is a convex functional. In view of Theorem \ref{p4-theorem005}, Theorem $\ref{Compctness}$ and Theorem $\ref{Close-to-convex-hull--Class-U}$,  we consider the function $f_0$ of the form
$$f_{0}(z)=\frac{z}{(1-xz)^{2}}=z+\sum\limits_{n=2}^{\infty}A_{n}z^{n},
$$ 
where $|x|=1$ and $A_{n}=nx^{n-1}$. Therefore, $\phi(f_{0})=A_{n}^{2}-A_{2n-1}=(n^{2}-2n+1)x^{2n-2}$ and hence $|\phi(f_{0})|=(n-1)^{2}$. For $x=e^{i\theta}$, a simple computation shows that
\begin{eqnarray*}
J(f_{0})&=&({\rm Re\,}(a_{n}))^{2}-{\rm Re\,}(a_{2n-1})\\
&=&n^{2}({\rm Re\,}(x^{n-1}))^{2}-(2n-1){\rm Re\,}(x^{2n-2})\\
&=&n^{2} \cos^{2}(n-1)\theta-(2n-1)\cos(2n-1)\theta\\
&=&n^{2} \cos^{2}(n-1)-(2n-1)(2\cos^{2}(n-1)\theta-1)\\
&=&\cos^{2}(n-1)\theta(n^{2}-4n+2)+2n-1\\
&\le & n^{2}-4n+2+2n-1 \qquad(\mbox{ since } n^{2}-4n+2>0 \mbox{ for } n\ge 4.) \\
&=& n^{2}-2n+1=(n-1)^{2}. 
\end{eqnarray*}
 Therefore,  $({\rm Re\,}(a_{n}))^{2}-{\rm Re\,}(a_{2n-1})\le (n-1)^{2}$ and hence 
${\rm Re\,}(a_{n}^{2}-a_{2n-1})\le (n-1)^{4}$. This implies that $|a_{n}^{2}-a_{2n-1}|\le(n-1)^{2}$ and the equality holds for the function 
of the form $f_{0}(z)=z/(1-xz)^{2}$, where $|x|=1$.  This completes the proof. 
\end{pf}

\section{Generalized Zalcman Conjecture for the class $\mathcal{U}$}\label{Section-4}

Let  $p(z)=1+\sum_{k=1}^{\infty}c_{k}z^{k} \in \mathcal{P}$ then from $(\ref{Abhi-Vasu-P1-Equ})$, it is easy to see that for all $n, m\in \mathbb{N}$,
$$\lambda c_{n-1}c_{m-1}-c_{n+m-1}=2\left(2\lambda\int\limits_{0}^{2\pi}e^{i(n-1)t}d\nu(t)\int\limits_{0}^{2\pi}e^{i(m-1)t}d\nu(t)-\int\limits_{0}^{2\pi}e^{i(n+m-2)t}d\nu(t)\right).
$$
By using Lemma \ref{Lemma-Ravichandran-2015},  we obtain
\begin{eqnarray}\label{Abhi-Vasu-P1-Equry}
&& \left|2\lambda\int\limits_{0}^{2\pi}e^{i(n-1)t}d\nu(t)\int\limits_{0}^{2\pi}e^{i(m-1)t}d\nu(t)-\int\limits_{0}^{2\pi}e^{i(n+m-2)t}d\nu(t)\right|\\[5mm]
&&\hspace{5cm} \le
\begin{cases}
 1,  & 0\le\lambda\le1\\[5mm]
 |2\lambda-1|,  & \mbox{ elsewhere }.\nonumber
\end{cases}
\end{eqnarray}

We now prove the generalized Zalcman conjecture for the class $\mathcal{U}$.
\begin{thm}\label{Generalized-Zalcman-Conjecture-Class-U} 
Let  $f\in \overline{co}\, \mathcal{U}$  be given by $(\ref{p4_i001})$. Then for $n,m \ge 2$
\[|a_{n}a_{m}-a_{n+m-1}|\le\begin{cases}
   n+m-1,&\mbox{ if }(n,m)\mbox{ is } (2,n), (m,2), (3,3), (3,4), (4,3)\\[5mm]
   (n-1)(m-1),&\mbox{ otherwise. }
   \end{cases}
   \]
   The second  inequality is sharp and the equality holds for the Koebe function and its rotations.
\end{thm}
\begin{pf}
%[{\bf Proof of Theorem   \ref{Generalized-Zalcman-Conjecture-Class-U}}] 
Let $f\in \overline{co}\,\mathcal{U}$ be given by $f(z)=z+\sum_{n=2}^{\infty}a_{n}z^{n}$. Then from Theorem $\ref{Close-to-convex-hull--Class-U}$,   
there exists a probability measure $\mu$ on $\partial\mathbb{D}$ such that 
$$f(z)=\int_{|x|=1}\frac{z}{(1-xz)^{2}}\, d\mu(x).
$$
Equivalently, there exists a probability measure $\nu$ on $[0,2\pi]$ such that 
$$f(z)=\int\limits_{0}^{2\pi}\frac{z}{(1-e^{it}z)^{2}}\, d\nu(t)$$ 
which can be written as, 
\begin{equation}\label{Abhi-Vasu-P1-Equ2}
z+\sum\limits_{n=2}^{\infty}a_{n}z^{n}=z+\sum\limits_{n=2}^{\infty}n\left(\int\limits_{0}^{2\pi}e^{i(n-1)t}d\nu(t)\right)z^{n}. 
\end{equation}
By comparing both the sides of $(\ref{Abhi-Vasu-P1-Equ2})$, we obtain
\begin{equation}\label{Abhi-Vasu-P1-Equ3}
a_{n}=n\int\limits_{0}^{2\pi}e^{i(n-1)t}\, d\nu(t).
\end{equation}
Using $(\ref{Abhi-Vasu-P1-Equ3})$, we obtain $$a_{n}a_{m}-a_{n+m-1}=nm\int\limits_{0}^{2\pi}e^{i(n-1)t}\, d\nu(t)\int\limits_{0}^{2\pi}e^{i(m-1)t}\, d\nu(t)-(n+m-1)\int\limits_{0}^{2\pi}e^{i(n+m-2)t}\, d\nu(t)$$
which can be written as
\begin{eqnarray}\label{Abhi-Vasu-P1-Equ4}
& a_{n}a_{m}-a_{n+m-1}=& n+m-1\left[\frac{nm}{n+m-1}\int\limits_{0}^{2\pi}e^{i(n-1)t}d\nu(t)\int\limits_{0}^{2\pi}e^{i(m-1)t}d\nu(t)\right.\\ \nonumber
&&\qquad -\left.\int\limits_{0}^{2\pi}e^{i(n+m-2)t}d\nu(t)\right].
\end{eqnarray}
Comparing $(\ref{Abhi-Vasu-P1-Equry})$ and $(\ref{Abhi-Vasu-P1-Equ4})$,  we obtain 
\begin{equation}\label{Abhi-Vasu-P1-Equ5}
\lambda=\frac{nm}{2n+2m-2}\, .
\end{equation}
For $\lambda\le1$, we see that $(\ref{Abhi-Vasu-P1-Equ5})$ implies $nm-2n-2m+2\le0$ which is equivalent to $(n-2)(m-2)\le2$.
If any of $n,m$ is $2$ and for the pairs $(n,m)=(3,3), (3,4), (4,3)$ the inequality $(n-2)(m-2)\le2$ holds. In all other choices of $n,m$ we get $\lambda>1$. Therefore, we have
 \begin{eqnarray}\label{Abhi-Vasu-P1-Equry1}
&& \left|\frac{nm}{n+m-1}\int\limits_{0}^{2\pi}e^{i(n-1)t}d\nu(t)\int\limits_{0}^{2\pi}e^{i(m-1)t}d\nu(t)-\int\limits_{0}^{2\pi}e^{i(n+m-2)t}d\nu(t)\right|\\[5mm]
 &&\hspace{5cm}  \le
 \begin{cases}
   1,  & (2,m), (n,2), (3,3), (3,4), (4,3) \\[5mm]
 % \displaystyle  
  \frac{(n-1)(m-1)}{n+m-1},  & \mbox{ other pairs of } (n,m). \nonumber
  \end{cases}
  \end{eqnarray}
   Therefore from $(\ref{Abhi-Vasu-P1-Equ4})$ and $(\ref{Abhi-Vasu-P1-Equry1})$,  we obtain
    \[|a_{n}a_{m}-a_{n+m-1}|\le\begin{cases}
   n+m-1,&\mbox{ if } (n,m) \mbox{ is } (2,n), (m,2), (3,3), (3,4), (4,3)\\[5mm]
   (n-1)(m-1),&\mbox{ otherwise.}
   \end{cases}
   \]
   The second inequality is sharp and equality holds for the Koebe function and its rotations. 
   \end{pf}
   
   \section{Zalcman Conjecture for certain subclass of Close-to-convex functions}\label{Section-5}
      In this section,  we consider the class of close-to-convex functions with argument $0$ and with respect to Koebe function $k(z)=z/(1-z)^{2}$. More precisely,  let 
      $$\mathcal{F}=\{f\in \mathcal{A}: {\rm Re\,}(1-z)^{2}f'(z)>0, \quad z\in \mathbb{D}\}.
      $$ 
      Clearly the functions in $\mathcal{F}$ are convex in the positive direction of the real axis. The region of variability of the class $\mathcal{F}$ has been studied  by Ponnusamy {\it et al} in 
      \cite{Samy-Vasu-Yanagihara}. In 2017, Ali and Vasudevarao \cite{Vasu-2017} obtained the sharp logarithmic  coefficient for functions in the class $\mathcal{F}$. \\ 
      
 Our aim is to solve the Zalcman conjecture  for $n=2$ and  generalized Zalcman Conjecture  for $n=2, m=3$.
      \begin{thm}\label{Abhi-Vasu-P1-subclass}
      Let $f\in\mathcal{F}$ given by $(\ref{p4_i001})$. Then
      \begin{enumerate}
      \item[(i)] $|a_{2}^{2}-a_{3}|\le 1$
      \item[(ii)] $|a_{2}a_{3}-a_{4}|\le 2$.
      \end{enumerate}
      These inequalities are sharp with equality for the Koebe function $k(z)=z/(1-z)^{2}$ and its rotations.
    \end{thm}
    \begin{pf}
    For $f\in\mathcal{F}$, let $g(z)=(1-z)^{2}f'(z)$. Clearly ${\rm Re\,}g(z)>0$. Since ${\rm Re\,}g(z)>0$ in $\mathbb{D}$, 
    there exists an analytic function $\phi:\mathbb{D}\to\mathbb{D}$  such that
     \begin{equation}\label{Abhi-Vasu}
    \phi(z)=\frac{g(z)-1}{g(z)+1}.
    \end{equation} 
    Clearly,  $\phi(0)=0$. Let
   $$\phi(z)=\sum\limits_{n=1}^{\infty}c_{n}z^{n}=c_{1}z+c_{2}z^{2}+c_{3}z^{3}+c_{4}z^{4}+\ldots.
   $$
   In view of the Schwarz Lemma, we get $|c_{1}|\le 1$. From $(\ref{p4_i001})$ we have
    $$g(z)=(1-z)^{2}f'(z)=(1+z^{2}-2z)(1+\sum\limits_{n=2}^{\infty}na_{n}z^{n-1}).$$
    From $(\ref{Abhi-Vasu})$, we have
    $$\phi(z)(g(z)+1)=g(z)-1.
    $$
    A simple computation shows that
    \begin{eqnarray}\label{Abhi-Vasu-00001}
    & a_{2}=& 1+c_{1}\\ \nonumber
   & a_{3}=& \frac{1}{3}\left(2c_{1}^{2}+4c_{1}+2c_{2}+3\right)\\ \nonumber
   &a_{4}=&\frac{1}{4}\left(2c_{1}^{3}+4c_{1}^{2}+4c_{1}c_{2}+6c_{1}+4c_{2}+2c_{3}+4\right). 
    \end{eqnarray} 

    Therefore using (\ref{Abhi-Vasu-00001}) we have
    \begin{eqnarray*}
    a_{2}^{2}-a_{3} &=& (1+c_{1})^{2}-\frac{1}{3}\left(2c_{1}^{2}+4c_{1}+2c_{2}+3\right)\\
    &=& \frac{1}{3}\left(3(1+2c_{1}+c_{1}^{2})-2c_{1}^{2}-4c_{1}-2c_{2}-3\right)\\
    &=& \frac{1}{3}\left(3+6c_{1}+3c_{1}^{2}-2c_{1}^{2}-4c_{1}-2c_{2}-3 \right)\\
    &=& \frac{1}{3}\left(c_{1}^{2}+2c_{1}-2c_{2}\right)
    \end{eqnarray*}
    which yields 
    $$|a_{2}^{2}-a_{3}|=\frac{1}{3}\left(|c_{1}^{2}+2c_{1}-2c_{2}|\right)\le \frac{1}{3}\left(|c_{1}|^{2}+2|c_{1}|+2|c_{2}|\right).$$
    It is well-known that $|c_{n}|\le 1-|c_{1}|^{2}$ for $n\ge 2$. Therefore,
    \begin{equation}\label{Abhi-Vasu-P1-equation1}
    |a_{2}^{2}-a_{3}|\le \frac{|c_{1}|^{2}+2|c_{1}|+2(1-|c_{1}|^{2})}{3}=\frac{2|c_{1}|+2-|c_{1}|^{2}}{3}.
    \end{equation}
    Let $x=|c_{1}|$ and $|c_{1}|\le 1$. Let $h(x)=2x+2-x^{2}$, where $0\le x\le 1$. Clearly, $h$ is non-negative in $[0,1]$ and $h'(x)=2-2x$ is also non-negative in $[0,1]$. 
    Therefore, $h$ is an increasing function and hence $h(1)=3$ is the maximum value. Hence $2|c_{1}|+2-|c_{1}|^{2}\le 3$. In view of this and  $(\ref{Abhi-Vasu-P1-equation1})$,  we obtain 
    $$|a_{2}^{2}-a_{3}|\le 1.
    $$
    A simple computation using  (\ref{Abhi-Vasu-00001})  gives
    \begin{eqnarray*}
    a_{2}a_{3}-a_{4}&=&
    \frac{1}{3}\left((1+c_{1})(2c_{1}^{2}+4c_{1}+2c_{2}+3\right)-\frac{1}{4}\left(2c_{1}^{3}+4c_{1}^{2}+4c_{1}c_{2}+6c_{1}+4c_{2}\right.\\
    &&\qquad\left.+2c_{3}+4\right)\\
    &=&\frac{1}{3}\left(2c_{1}^{3}+6c_{1}^{2}+2c_{1}c_{2}+7c_{1}+2c_{2}+3\right)-\frac{1}{4}\left(2c_{1}^{3}+4c_{1}^{2}+4c_{1}c_{2}+6c_{1}+\right.\\
    &&\qquad\left. 4c_{2}+2c_{3}+4\right)\\
    &=&\frac{1}{12}\left(2c_{1}^{3}+12c_{1}^{2}-4c_{1}c_{2}+10c_{1}-4c_{2}-6c_{3}\right)\\[5mm]
    &=&\frac{1}{6}\left(c_{1}^{3}+6c_{1}^{2}-2c_{1}c_{2}+5c_{1}-2c_{2}-3c_{3}\right).
    \end{eqnarray*}
    Therefore,
    \begin{eqnarray*}
    |a_{2}a_{3}-a_{4}| &\le& \frac{1}{6}\left(|c_{1}|^{3}+6|c_{1}|^{2}+2|c_{1}||c_{2}|+5|c_{1}|+2|c_{2}|+3|c_{3}|\right)\\
    &\le& \frac{1}{6}\left(|c_{1}|^{3}+6|c_{1}|^{2}+2|c_{1}|(1-|c_{1}|^{2})+5|c_{1}|+2(1-|c_{1}|^{2})+3(1-|c_{1}|^{2})\right)\\
    &=& \frac{1}{6}\left(-|c_{1}|^{3}+|c_{1}|^{2}+7|c_{1}|+5\right).
    \end{eqnarray*}
    Thus we obtain
    \begin{equation}\label{Abhi-vasu-P1-equation2}
    |a_{2}a_{3}-a_{4}|\le \frac{1}{6}\left(-|c_{1}|^{3}+|c_{1}|^{2}+7|c_{1}|+5\right).
    \end{equation}
     Let $x=|c_{1}|$, $|c_{1}|\le 1$ and $L(x)=-x^{3}+x^{2}+7x+5$, where $0\le x\le 1$. Then $L'(x)=-3x^{2}+2x+7$. It is easy to see that both the functions 
     $L$ and $L'$ are non-negative in $[0,1]$. Therefore,  $L$ is increasing in $[0,1]$ and $L(1)=12$ is the maximum value. 
     Hence  $-|c_{1}|^{3}+|c_{1}|^{2}+7|c_{1}|+5\le 12$. Using $(\ref{Abhi-vasu-P1-equation2})$,  we obtain 
    $$|a_{2}a_{3}-a_{4}|\le 2.$$
    This completes the proof.
 \end{pf}
 
%%%%%%%%%%%%%%%%%%%%
%%%%%%%%%%%%%%%%%%%%%
%%%%%%%%%%%%%%%%%%%%%

\vspace{4mm}
\noindent\textbf{Acknowledgement:} 
The first author thank SERB-MATRICS 
and the second author thank PMRF-MHRD, Govt. of India for their support.

\end{document}